\documentclass[12pt]{article}
\usepackage[latin1]{inputenc}
\usepackage{mathptmx}
\usepackage{a4wide}
\usepackage{amsmath}
  \usepackage{paralist}
  \usepackage{graphicx}
\usepackage{afterpage}
 \usepackage[colorlinks=true]{hyperref}

\title{An iterative method for the canard explosion\\ in general planar systems}
\author{Morten Br{\o}ns\\
Department of Mathematics\\ Technical University of Denmark\\2800
Lyngby, Denmark\\ \texttt{m.brons@mat.dtu.dk}\vspace*{2em}}

\date{Paper presented at the 9th AIMS Conference on Dynamical Systems,
Differential Equations and Applications, Orlando, Florida, USA
July 1 - 5, 2012\\[1em] Manuscript date \today}

\begin{document}
\maketitle

\begin{abstract}
  The canard explosion is the change of amplitude and period of a
  limit cycle born in a Hopf bifurcation in a very narrow parameter
  interval. The phenomenon is well understood in singular perturbation
  problems where a small parameter controls the slow/fast
  dynamics. However, canard explosions are also observed in systems
  where no such parameter is present. Here we show how the iterative
  method of Roussel and Fraser, devised to construct regular slow
  manifolds, can be used to determine a canard point in a general
  planar system of nonlinear ODEs. We demonstrate the method on the
  van der Pol equation, showing that the asymptotics of the method is
  correct, and on a templator model for a self-replicating system.
\end{abstract}
\clearpage

\section{Introduction}
Since the original discovery of canards in the van der Pol equation
more than 30 years ago \cite{benoit-callot1981:chass}, they have been
identified in numerous systems of nonlinear ODEs. A \emph{canard} is a
trajectory which stays close to a repelling slow manifold for an
extended amount of time. Canards play a key role as parts of
transitional limit cycles linking small cycles born in a Hopf
bifurcation with large relaxation oscillations when a parameter is
varied. Since this transition typically takes place over a very short
parameter interval, and easily may be mistaken for a discontinuous
event, the phenomenon has been denoted a \emph{canard explosion}. 

The mathematical theory for canards is well-established for singular
perturbation systems of the form
\begin{equation}
    \label{eq:1}
    \dot x = f(x,y,c,\epsilon), \quad \dot y = \epsilon g(x,y,c,\epsilon),
\end{equation}
where $\epsilon$ is a small parameter and $c$ is a bifurcation
parameter see e.g.\ \cite{benoit-callot1981:chass,
  eckhaus1983:relax-oscil-includ,
  krupa-szmolyan2001:relax-canar-explos}. In particular, asymptotic
expansions in terms of $\epsilon$ of the canard point $c_c$, the
parameter values where the longest canards exist, can be obtained
\cite{zvonkin-shubin1984:non-analy-singul,
  brons2005:relax-oscil-canar}. However, canard explosions have also
been observed in many systems that do not have an explicit slow/fast
structure with a well-defined small parameter $\epsilon$,
\begin{equation}
  \label{eq:2}
    \dot x = F(x,y,c), \quad \dot y = G(x,y,c).
\end{equation}
In some cases a small parameter can be identified after coordinate
transformations \cite{brons2011:canar-explos-limit}, while in other
cases an artificial parameter must be introduced to allow an
asymptotic expansion \cite{brons-bar-eli1994:asymp-analy,
  brons-sturis2001:explos-limit-cycles-chaot-waves}. After the
expansion, the artificial parameter is set to one to recover the
original system.

While these approaches have been successful, they are somewhat ad-hoc,
and it would be of interest to establish a systematic approach to
identify and locate canard explosions in general systems of the form
Eqns.~\eqref{eq:2}. The purpose of the present paper is to provide
such a procedure. It is a simple modification of the iterative method
by Fraser and Roussel \cite{fraser1988:stead-state-equil-approx,
  roussel-fraser1990:geomet-stead-state-approx} for finding slow
manifolds. We show that for the van der Pol equation with a
distinguished small parameter the method gives the correct asymptotic
result. For the templator model \cite{brons2011:canar-explos-limit}
with no small parameter we get an excellent agreement between the
canard point found from simulations and the lowest-order canard point
from the method.

\begin{figure}[th]
  \centering
  \includegraphics[width=\textwidth]{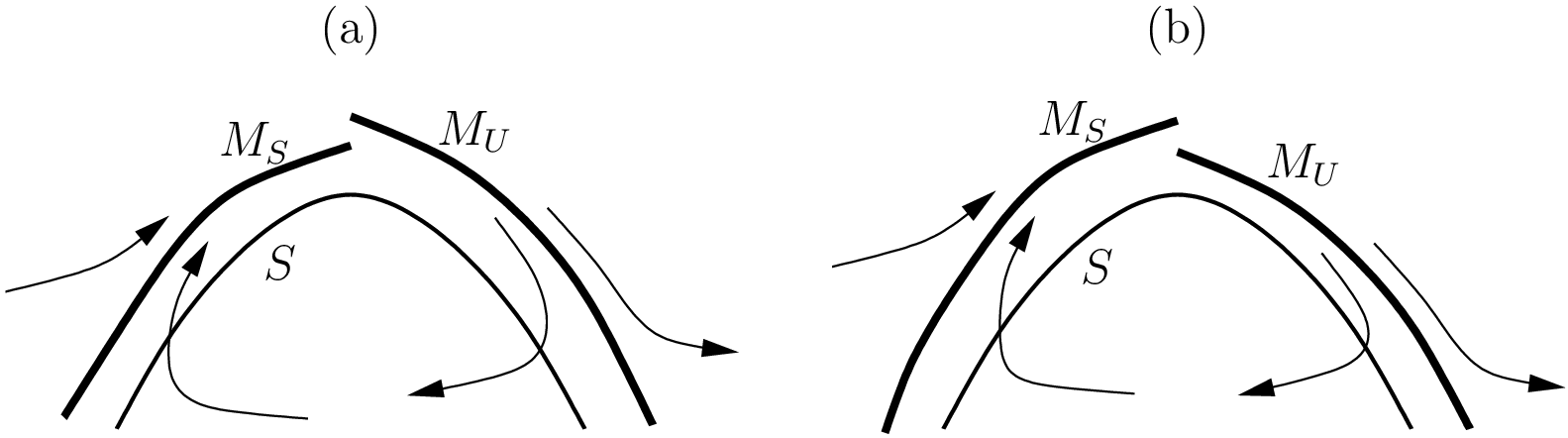}
  \caption{Canard explosion in a singular perturbation system,
    Eqns.~\eqref{fig:1}. As trajectories cross the fold of the critical
  manifold $S$, they are either repelled down or up, depending on the
  relative position of the slow manifolds $M_S$ and $M_U$.}
  \label{fig:1}
\end{figure}

\section{The canard explosion}
\label{sec:canard-explosion}
Here we briefly review the basics of the theory for the canard
explosion for Eqns.~\eqref{eq:1}. The curve defined by $f(x,y,c,0)=0$
is denoted \emph{the critical manifold} $S$. For $\epsilon=0$, $S$
consists of fixed points and assuming that it has a fold, the local
phase portrait is as shown in Fig.~\ref{fig:1}. For $\epsilon>0$ it
follows from standard Fenichel theory (see e.g.\
\cite{verhulst2005:method-applic-singul-pertur}) that on the stable
side of $S$ an attracting slow manifold $M_S$ exists and on the
unstable side a repelling slow manifold $M_U$ exists. The existence
and the normal hyperbolicity of these manifolds is guaranteed by the
theory away from the fold point only. However, as trajectories they
may be extended across the fold point. In general, $M_S$ and $M_U$
will be distinct, but for a special value of $c=c_c$ they may coincide
and form a single trajectory, a canard. Clearly, the shape of a limit
cycle will change dramatically if the parameter is varied across
$c_c$. If $M_U$ is above $M_S$ as in Fig.~\ref{fig:1}(a) only small
limit cycles will be possible. If $M_U$ is below $M_S$ as in
Fig.~\ref{fig:1}(b) the limit cycles will be large.

The single trajectory $M_S=M_U$ and the corresponding parameter value
$c_c$ can be found asymptotically. For the equation for the trajectories
\begin{equation}
  \label{eq:3}
  f(x,y,c,\epsilon)\frac{dy}{dx} = \epsilon g(x,y,c,\epsilon)
\end{equation}
a Poincaré-Lindstedt series is inserted,
\begin{equation}
  \label{eq:4}
  y = y_0 + y_1\epsilon + y_2\epsilon^2 + \cdots, \quad c_c = c_0 +
  c_1\epsilon  + c_2\epsilon^2 + \cdots.
\end{equation}
Collecting terms of the same order in $\epsilon$ algebraic equations
for the $y_k$ are obtained. These will in general have a singularity
at the fold point, but there will be a choice of $c_{k-1}$ such that
this singularity cancels and $y_k$ is well-defined at the fold point. 
This choice defines the canard point and $y_k$ is the corresponding
canard solution.

\section{A general iterative procedure}
\label{sec:gener-iter-proc}
For Eqns.~\eqref{eq:2} we can also write down the equation for the
trajectories,
\begin{equation}
  \label{eq:5}
  F(x,y,c)\frac{dy}{dx} = G(x,y,c).
\end{equation}
Following Fraser and Roussel
\cite{fraser1988:stead-state-equil-approx,
  roussel-fraser1990:geomet-stead-state-approx}, we solve this
equation for $y$ algebraically,
\begin{equation}
  \label{eq:6}
  y = \Phi\left(x,\frac{dy}{dx},c\right).
\end{equation}
Clearly, it must be assumed that such a solution exists, at least
locally. From this equation an iterative procedure can be established,
\begin{equation}
  \label{eq:7}
  y_k = \Phi\left(x,\frac{dy_{k-1}}{dx},c\right).
\end{equation}
To start the iteration we choose $y_0$ such that $F(x,y_0(x,c),c)=0$,
that is, the $\infty$-isocline. Again, we must assume that this
equation can be solved for $y_0$. Other choices will be possible, but
we do not have space here to discuss this issue.  Typically $y_k$ will
have a singularity. Since $y_k$ depends on $c$, we will choose the
value in each step such that this singularity cancels. This defines
the procedure for finding canards and canard points for
Eqns.~\eqref{eq:2}.

\section{The van der Pol equation}
\label{sec:van-der-pol}
We now demonstrate the iterative method on the van der Pol equation
\begin{equation}
  \label{eq:8}
  \dot x = y - (x^3/3-x), \quad \dot y = \epsilon(c-x).
\end{equation}
This system has a canard explosion for $c$ close to 1 when $\epsilon$
is small and positive. The procedure from
\S~\ref{sec:canard-explosion} yields for the canard point
\cite{zvonkin-shubin1984:non-analy-singul} 
\begin{equation}
  \label{eq:9}
  c_c = 1 - \frac{1}{8}\epsilon - \frac{3}{32}\epsilon^2
  -\frac{173}{1024}\epsilon^3 + {\mathcal O}(\epsilon^4).
\end{equation}

The iterative procedure Eqn.~\eqref{eq:7} is defined by
\begin{equation}
  \label{eq:10}
  y_{k+1}=  x^3/3 - x + \epsilon \frac{c-x}{y'_k}
\end{equation}
with starting point
\begin{equation}
  \label{eq:11}
  y_0 = x^3/3 - x.
\end{equation}

\subsection{A numerical example}
\label{sec:vdpnumex}
We consider first the van der Pol system Eqns.~\eqref{eq:8} with
$\epsilon=0.1$. The asymptotic formula for the canard point
Eqn.~\eqref{eq:9} yields $c_c = 0.986394$.

The iterative process runs as follows: From Eqn.~\eqref{eq:10} we find
\begin{equation}
\label{eq:18}
y_1 = x^3/3 - x + \frac{1}{10}\frac{c-x}{x^2-1}.
\end{equation}
This has a singularity  at $x=1$ (and also at $x=-1$, but this is not
of interest here) which is removed by choosing $c=1$, which, then, is
the first approximation to the canard point. The relative deviation
from the asymptotic value is $1.38\%$. With this choice of $c$ we have 
\begin{equation}
\label{eq:19}
y_1 = x^3/3 - x - \frac{1}{10}\frac{1}{x+1}
\end{equation}
and a further iteration yields 
\begin{equation}
\label{eq:20}
y_2 = x^3/3 - x + \frac{(x+1)^2(c-x)}{p_2}
\end{equation}
where 
\begin{equation}
\label{eq:21}
p_2 = 10x^4 + 20x^3 - 20x - 9.
\end{equation}
The polynomial $p_2$ has two real roots, $x_1=-0.603433$ and
$x_2=0.987258$. We remove the singularity of $y_2$ at the latter point
by choosing $c=x_2=0.987258$, which is the second approximation to the
canard point. This deviates from the asymptotic value by $0.09\%$. 

By factorization we get 
\begin{equation}
\label{eq:22}
p_2 = (x-x_2)q_2
\end{equation}
where 
\begin{equation}
\label{eq:23}
q_2 = 10x^3 + 29.8726x^2 + 29.4919x + 9.11616
\end{equation}
such that 
\begin{equation}
\label{eq:24}
y_2 = x^3/3 - x - \frac{(x+1)^2}{q_2}.
\end{equation}
By iteration we find $y_3$, which is a rational function where the
denominator is a polynomial of degree 8 in $x$ but independent of $c$. The
real roots are $-1.24503$, $-0.999999$, $-0.389117$ and
$x_3=0.986481$. The numerator of $y_3$ is a polynomial of degree 11 in
$x$ but linear in $c$. The singularity at $x_3$ can be canceled by
choosing $c=x_3=0.986481$. This gives yet an improvement of the canard
point, as the deviation from the asymptotic value is now down to
$0.009\%$. Clearly, the procedure can be continued any number of
times.

\subsection{Asymptotic analysis}
\label{sec:asymp}
The structure of the van der Pol equation is sufficiently simple to
allow an asymptotic analysis in the limit $\epsilon\rightarrow 0$ of
the iterative procedure. For a general $\epsilon$, we get in the first
iteration 
\begin{equation}
  \label{eq:12}
  y_1 =  x^3/3 - x + \frac{c-x}{x^2-1}\epsilon.
\end{equation}
As before, we eliminate the singularity at $x=1$ by choosing $c=1$
such that
\begin{equation}
  \label{eq:13}
  y_1 =  x^3/3 - x-\frac{1}{x+1}\epsilon.
\end{equation}
The next iteration gives
\begin{equation}
  \label{eq:14}
  y_2 =  x^3/3 - x + \frac{(x+1)^2(c-x)}{p_2}\epsilon
\end{equation}
where \begin{equation}
\label{eq:25}
p_2 = x^4+2x^3-2x-1+\epsilon.
\end{equation}
The function $y_2$ has a singularity at $x=x_2$ where $x_2$ is a root
of $p_2$ and the singularity cancels if $c=c_2=x_2$. The polynomial
$p_2$ has two roots for $\epsilon<27/16$, so the construction of the
canard trajectory only works under this condition. When
$\epsilon=27/16$, $p_2$ has a double root at $x=1/2$. We choose as
$x_2$ the root which is greater than $1/2$. By a Taylor expansion, one
easily finds 
\begin{equation}
  \label{eq:16}
  x_2 = c_2 = 1-\frac{1}{8}\epsilon-\frac{3}{128}\epsilon^2 -\frac{15}{2048}\epsilon^3 + O(\epsilon^4).
\end{equation}
This agrees with \eqref{eq:9} to $O(\epsilon)$, but not to
$O(\epsilon^2)$.  Inserting $c=c_2$ from Eq.~\eqref{eq:16} in $y_2$
and iterating in Eq.~\eqref{eq:10} yields $y_3$ as a rational
function. In this, we let $c=c_0+c_1\epsilon+c_2\epsilon^2+\cdots$ and
in a Taylor expansion the first terms are 
\begin{equation}
\label{eq:15}
y_3 = \frac{c_0-x}{(x-1)(x+1)}\epsilon +
\frac{c_1x^4+2c_1x^3-(2c_1+1)x-c_0-c_1}{(x-1)^2(x+1)^4}\epsilon^2 + \cdots 
\end{equation}
By choosing $c_0=1$ and $c_1=-1/8$ the singularities at $x=1$ in the
first two terms cancel. Proceeding to the term of order $\epsilon^3$
(we omit the rather long expression) one cancels a singularity by
choosing $c_2=-3/32$. Continuing this way, we find an approximation to
the canard point as
\begin{equation}
\label{eq:17}
c = 1-\frac{1}{8}\epsilon-\frac{3}{32}\epsilon^2 -\frac{75}{1024}\epsilon^3 + O(\epsilon^4)
\end{equation}
This agrees with \eqref{eq:9} to $O(\epsilon^2)$, but not to
$O(\epsilon^3)$. Again, we may continue this procedure to any order,
in each step correcting a term in the asymptotic expansion of the
canard point.

\afterpage{\clearpage}

\section{The templator}
\label{sec:templator}
The templator is a mathematical model for the kinetics of a
self-replicating chemical system. The reactions are
\begin{gather}
  X_0 \rightarrow X \notag \\
  X + X \rightarrow T \notag \\
  X + X + T \rightarrow T + T \notag \\
  T \rightarrow P \notag
\end{gather}
The key process the third one where a dimer $T$ acts as a templates
and catalyzes its own production from a monomer $X$. In dimensionless
variables the model can be written

\begin{figure}[thb]
  \centering
  \begin{tabular}{cc}
    (a) & (b)\\
\includegraphics[width=0.45\textwidth]{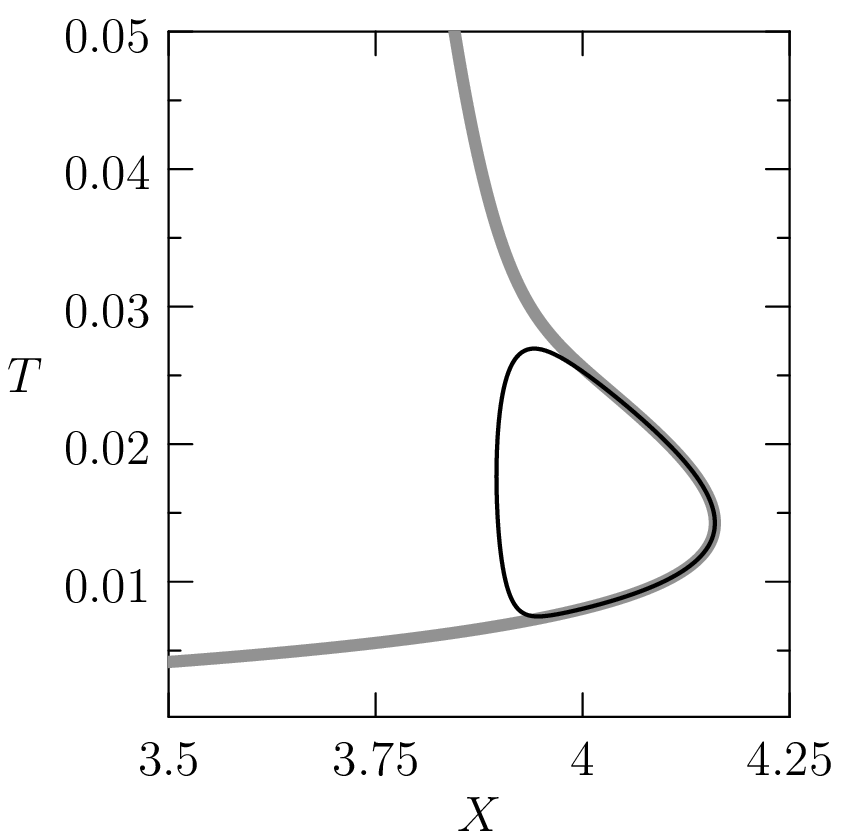} &
\includegraphics[width=0.45\textwidth]{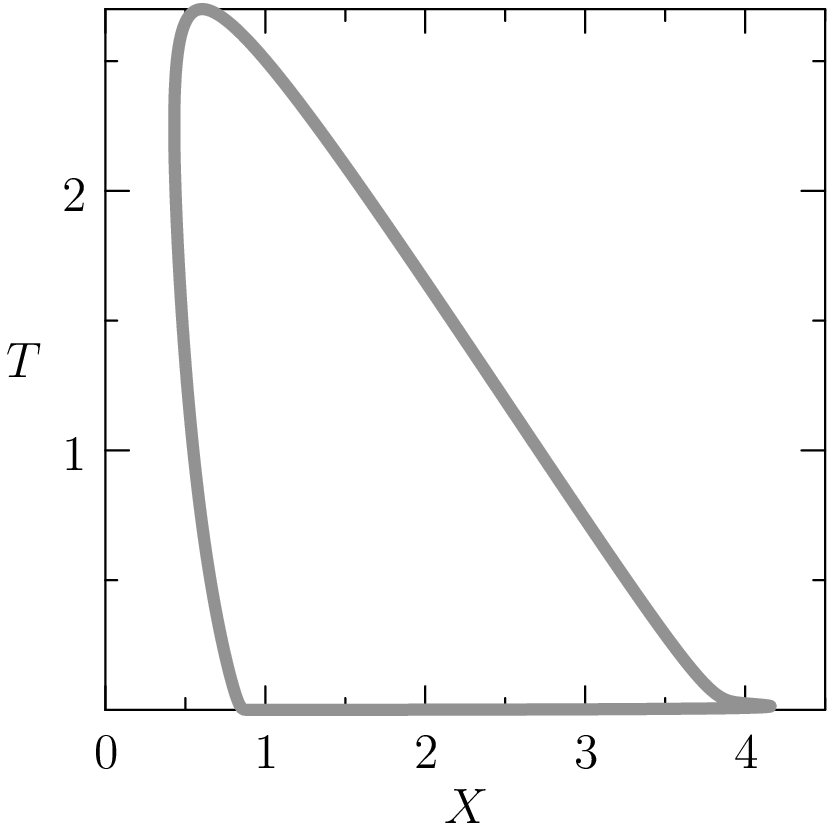}
\end{tabular}
\begin{tabular}{c}
 (c) \\
\includegraphics[width=0.45\textwidth]{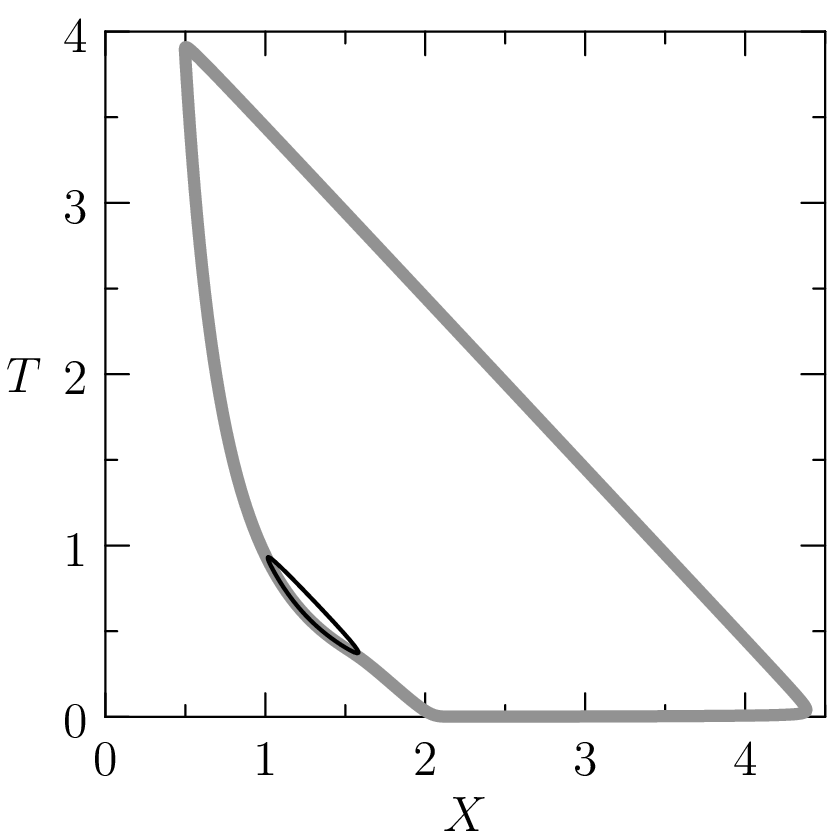} 
  \end{tabular}
  \caption{Simulations of the templator model Eqns.~\eqref{eq:26}. (a)
    The black curve is the limit cycle for $r=0.419940$, the gray
    curve is is a part of the limit cycle for $r=0.419945$.  (b) The
    full large limit cycle from panel (a). Note the differences in the
    scales. (c) The gray curve is the limit cycle for $r=0.96755$, the
    black curve the limit cycle for $r=0.96756$.}
  \label{fig:tamplator}
\end{figure}

\begin{subequations}
  \label{eq:26}
  \begin{align}
    \frac{dX}{dt} &= r - k_u X^2 - k_T X^2T,  \\
    \frac{dT}{dt} &= k_u X^2 + k_T X^2T - \frac{qT}{K+T}.
  \end{align}
\end{subequations}
The last step in the reaction is modeled as an enzymatic reaction with
Michaelis-Menten kinetics. For further details on the model and its
biological significance see
\cite{beutel-peacock-lopez2007:compl-dynam-cross-self-mechan,
brons2011:canar-explos-limit} and references therein.

In the following we fix the parameters $k_u=0.01$, $k_T=1$, $q=1$,
$K=0.02$ and consider $r$ as a bifurcation parameter. In
\cite{beutel-peacock-lopez2007:compl-dynam-cross-self-mechan,
  brons2011:canar-explos-limit} it is shown numerically that two
canard explosions occur in the model. One is at $r=0.419942$ where a
small limit cycle explodes as $r$ increases. The large limit cycle
persists until $r=0.967555$, where it turns into a small cycle in
another canard explosion. See Fig.~\ref{fig:tamplator}. There is no
obvious small parameter in the equations, so the standard asymptotic
approach for Eqns. \eqref{eq:1} does not work. However, in
\cite{brons2011:canar-explos-limit} it is shown that it is possible to
account for the two canard explosions by two different scalings of the
equations. Here we show that the iterative method described in this
paper can be applied directly on the unscaled equations.

The equation for the trajectories is 
\begin{equation}
\label{eq:27}
\left(k_uX^2 + k_TX^2T - \frac{qT}{K+T}\right)\frac{dX}{dT} = r -
k_uX^2 - k_TX^2T.
\end{equation}
This is a quadratic equation in $X$. Choosing the positive solution,
we get for the iteration process
\begin{equation}
\label{eq:28}
X_{k+1} = \sqrt{\frac{r(K+T)+ qTX_k'}{(X_k'+1)(k_u+k_TT)(K+T)}}.
\end{equation}
To start the iteration, we choose the $dT/dt=0$ isocline as the
initial approximation, 
\begin{equation}
\label{eq:29}
X_0 = \sqrt{\frac{qT}{(k_u+k_TT)(K+T)}}.
\end{equation}
The expression for $X_1$ which is quite complicated has a denominator
which is independent of $r$. It has two zeroes, $T=0.0143454$ and
$T=0.599393$. Inserting these in the numerator of $X_1$ and requiring
that it is zero to cancel the singularity yields a linear equation to
determine $r$ with solutions $r=0.417681$ and $r=0.967710$
respectively. The first canard point deviates from the numerically
determined one by $0.6\%$, while the latter deviates with
$0.02\%$. Hence, a very accurate determination of the canard point is
achieved in the very first iteration.

\section{Conclusions}
\label{sec:conclusions}
We have demonstrated that a very simple iteration procedure can be
used to determine canard points in general planar dynamical systems
with no distinguished small parameter. We have shown for the van der
Pol equation that we obtain an asymptotically correct result in the
limit of $\epsilon\rightarrow 0$, and we conjecture this is a
general result for problems on the classical singular perturbation
form. For the more complex templator system, the method successfully
found the two canard points in one iteration. In the  analytical
approach \cite{brons2011:canar-explos-limit} different scalings were
needed to find the two canard points.

It is interesting to note that for the van der Pol equation an upper
bound for the small parameter $\epsilon$ for canard explosion to occur
was found. Recently a bound of $\epsilon< 1/4$ was found from
consideration of the curvature of the trajectories
\cite{desroches-jeffrey2011:canar-curvat}. The present bound is more
conservative, and it would be interesting to obtain a clearer
understanding of the relation of the two approaches to canard
explosion. Furthermore, the present iterative procedure provides a new
view on canard explosion which may lead to a more general
understanding on the specific conditions needed for a planar dynamical
system without a small parameter to experience a canard
explosion. Since systems with canard explosions of this kind are
abundant in the applications this seems to be of fundamental
interest. Work along these lines is in progress, and will be reported
elsewhere.

\bibliography{mb,canard,mybooks}
\bibliographystyle{plain}

\end{document}